\documentclass[12pt, leqno]{article}
\usepackage{latexsym}

\newtheorem{cor}{Corollary}[section]
\newtheorem{thm}[cor]{Theorem}
\newtheorem{lem}[cor]{Lemma}
\newtheorem{rem}[cor]{Remark}
\newtheorem{defi}[cor]{Definition}

\newtheorem{exa}[cor]{Example}

\newtheorem{str}[cor]{Strategy}

\newcommand{\RE}{\rm I\mskip-4mu R}

\newcommand{\lla}{ ( }
\newcommand{\rra}{ ) }

\newcommand{\M}{{\cal M}}

\newcommand{\N}{{\cal N}_A}

\newcommand{\U}{{\cal U}(A)}
\newcommand{\UT}{{\cal L}(A)}

\newcommand{\G}{{\cal G}}

\newcommand{\at}{T}

\newcommand{\af}{F}

\newcommand{\msn}{\medskip \noindent}

\begin{document}

\bibliographystyle{plain}

\title{Block LU Factorizations of M--matrices}

\author{J. J. McDonald \thanks{Work supported by an NSERC Research Grant}
\\Department of Mathematics and Statistics
\\University of Regina
\\Regina, Saskatchewan, Canada
\\S4S 0A2
\and H. Schneider \thanks {Work supported by NSF Grant
           DMS-9123318 and DMS-9424346.}
\\Department of Mathematics
\\University of Wisconsin
\\ Madison, Wisconsin 53706}
\date{19  AUGUST 1997}
\maketitle

\newpage
\begin{abstract}
It is well known that any nonsingular M--matrix admits an LU factorization
into M--matrices (with $L$ and $U$ lower and upper triangular respectively)
and any singular M--matrix is permutation similar to an M--matrix which admits
an LU factorization into M--matrices. Varga and Cai establish necessary and
sufficient conditions for a singular M--matrix (without permutation)
to allow an LU factorization with $L$ nonsingular.
We generalize these results in two directions. First, we find necessary
and sufficient conditions for the existence of an LU factorization of a
singular M-matrix where $L$ and $U$ are both permitted to be singular.
Second, we establish the minimal block structure that a block LU
factorization of a singular M--matrix can have when $L$ and $U$
are M--matrices.

\end{abstract}

\section{Introduction}

It was shown by Fiedler-Ptak, \cite{FP}, that any nonsingular
M--matrix $A$ admits an LU factorization, $A=LU$, where
$L$ is a nonsingular lower triangular M--matrix, and
$U$ is a nonsingular upper triangular M--matrix. Kuo,
\cite{K}, proved that any (singular) irreducible
M--matrix $A$ admits an LU factorization, $A=LU$, where
$L$ is a nonsingular lower triangular M--matrix, and
$U$ is a upper triangular M--matrix, and she gave an example to show
that not every singular M--matrix admits an $LU$ factorization of this type.
In \cite{VC}, Varga and Cai establish necessary and sufficient
conditions in terms of the directed graph $\G(A)$ of $A$
for a singular M--matrix to allow
a LU factorization into M--matrices with
$L$ nonsingular.

In this paper we consider the case where the
conditions outlined in \cite{VC} may not be satisfied.
We generalize these results in two directions. First, we find necessary
and sufficient conditions for the existence of an LU factorization of a
singular M--matrix where both $L$ and $U$ are permitted to be singular.
Second, we wish to factor an M--matrix $A$ into $A=LU$, where $L$ and
$U$ are M--matrices which are as close to lower and upper triangular as
possible. Our goal is to minimize the number of nonzeros above the diagonal
of $L$ and below the diagonal of $U$ in a factorization $A = LU$ and/or
to optimize on their placement. Our approach to this problem is to
minimize the appropriate access relationships in the digraphs $\G(L)$
and $\G(U)$ of $L$ and $U$ respectively.
In order to establish the
minimal block structure that a block LU factorization of a singular
M--matrix can have we need to give careful
definitions of what is meant by the {\em
block lower (upper) triangular self--partition} of a matrix.
These partitions are minimal in the set of partitions that lead to
a block lower (upper) triangular matrix
 and they are solely determined by the zero/non-zero pattern of the
matrix and thus do not depend on some assumed prior partitioning of the
matrix.
We use the term {\em block
factorization} to indicate that we are interested in the block
structure (without permutation) of the factors involved.

We now describe our paper in more detail.
Our definitions are contained in Section 2.

In Section 3 we
examine both LU and block LU factorizations of M--matrices
where $L$ and $U$ are permitted to be singular. Thus, in
 Example \ref{eg} we
provide an M--matrix $A$ which has an LU factorization
into M--matrices
only when $L$ and $U$ are both singular.  In Theorem
\ref{sbs} we identify the minimum access relationships which
must be present in $\G(L)$ and $\G(U)$
and thus also identify the
minimum sizes of the diagonal blocks of $L$ and $U$ in a block LU
factorization.
In Theorem \ref{ebs} we show that these minimum block structures
can be achieved.  Returning
to (elementwise) LU factorizations
in Theorem \ref{eq2}, we there characterize
the M--matrices $A$ which
admit an LU factorization into (possibly
singular) M--matrices. We give two strategies, one for finding
a desirable block LU factorization, and the other for choosing
a permutation matrix $P$ so that $PAP^T$ has an LU factorization.

In Section 4
we consider the
case where $L$ is a nonsingular M--matrix and
$U$ is an M--matrix.
Using definitions and a result
from McDonald \cite{M}, we actually examine a
slightly larger class - factorizations of an M--matrix
$A$ for
which $L$ is nonsingular and inverse nonnegative,
$L^{-1}$ is class nonsingular for $A$,
and $U$ is an
M--matrix.  In Theorem \ref{dtf} we show that for such
a factorization,  certain access
relationships must be present in $\G(U)$ and thus put a
lower limit on the number of nonzeros
below the diagonal of $U$.
In Theorem \ref{tt} we
show that we can attain these minimum access relationships
using a factorization for which $L$ is
a lower triangular nonsingular M--matrix.

In Section 5
we examine factorizations of $A$ into $LBU$ where
$L$ is a nonsingular lower triangular M--matrix,
$U$ is a nonsingular upper triangular M--matrix,
and $B$ is a block diagonal M--matrix.
Since $L$ and $U$
are nonsingular lower and upper triangular M--matrices,
standard methods could be used to solve the parts of the
system in which they are involved.
Since $B$ is generally a sparse matrix, it may be possible to
use specialized techniques in this area to solve the part of
the system involving $B$.

We
remark that analogous results hold
for UL--factorizations and that the techniques used in this
paper can also be used to identify the minimum structures
of $U$ and $L$ in this instance.

Our results are of a graph theoretical nature, but we express the hope
that our identification of possible block LU and LBU factorizations
of singular M-matrices will extend applications and numerical
implementations of block LU and LBU factorizations such as those
discussed in \cite[Chapter 12]{H}.

\section{Definitions}
We begin with some standard definitions. Let $n$ be a fixed positive
integer. We write $ \langle n \rangle=\{1,\ldots ,n \}.$

Throughout Section 2, $X=[x_{ij}]$ will denote a matrix in
$\RE ^{nn}$.

We say $X$ is:

\noindent
{\it positive} $(X \gg 0)$ if
$x_{ij}>0, $ for all $ i,j \in \langle n
\rangle$;

\noindent {\it semipositive}
$(X>0)$ if $x_{ij} \geq 0, $ for all $ i,j
\in \langle n \rangle$ and $X \ne 0$; and

\noindent {\it nonnegative}
$(X \geq 0)$ if $x_{ij} \geq 0, $ for all $
i,j \in \langle n \rangle$.

We say $X$  is a {\it Z--matrix} if
$X= \alpha I-P$ for some $\alpha  \in
\RE$ with $P$ nonnegative.  If
in addition, $\alpha $ is greater than or equal to the
spectral radius of $P$, then we say
$X$ is an {\it M--matrix}.  We denote the class of
$n\times n$ M--matrices by
$ \M$.

For any $J\subseteq\langle n\rangle$, we let
$$max(J)=max\{j\in J\},$$
$$ min(J)=min\{j\in J\},$$
$$J^{+}=\{l\in\langle n\rangle \  | \ l>max(J)\},$$
$$J^{-}=\{l\in\langle n\rangle \  | \ l<min(J)\},$$
$$J'=\{l\in\langle n\rangle \  | \ l\not\in J\}.$$
$$|J|=\mbox{ number of elements in } J.$$

For any $J,K \subseteq \langle n \rangle$ , we
write
$X_{JK}$ to represent the
submatrix of $X$ whose rows are indexed by
the elements of
$J$ and whose columns are
indexed by the elements of $K$, where the elements of
$J$ and $K$ are arranged in ascending order.

We call  the pair $\Gamma = (V,E)$ a {\it directed graph},
where $V$ is a finite set, and $E\subseteq V\times V$.
A {\it path} from $j$
to $k$ in $\Gamma $ is a sequence of vertices $
j=r_1,r_2,...,r_t=k$, with
$(r_i, r_{i+1}) \in E$, for $i=1,...,t-1$.   A path for which the
vertices are pairwise distinct
is called a {\it simple path}.  The
empty path is considered to be
a simple path linking every vertex to
itself. If there is a
path from $j$ to $k$, we say that $j$ has {\it
access} to $k$.   If $j$ has access to $k$ and $k$ has
access to $j$, we say $j$ and $k$ {\it communicate}.
The communication relation is an equivalence relation, and
hence we may partition V into equivalence classes, which
we will refer to as the {\it classes} of $\Gamma$. We define the
closure  $\Gamma$ by $\overline {\Gamma}= (V,F)$, where
$V= \langle n \rangle $ and $F = \{ (i,j) \  \vert \  i $ has
access to $ j $
in $\Gamma \} $.

We define the {\it (directed) graph} of $X$ by $\G (X)=(V,E)$,
where $V= \langle n \rangle $
and $E= \{ (i,j) \  \vert \  x_{ij} \ne 0 \} $.
It is well known that the
classes of $\G (X)$ correspond to the irreducible
components of $X$.  For any class $J$ of $G(X)$
we say that $J$ is a {\it
singular class} of $X$ if $X_{JJ}$ is singular,
and we
say that $J$
is a {\it nonsingular class }of $X$ if $X_{JJ}$ is
nonsingular.

We commonly view a matrix $X$ as a
block matrix $$X=\left[\begin{array}{ccc}
X_{11}&\ldots&X_{1p} \\
\vdots&\ddots&\vdots         \\
X_{q1}&\ldots&X_{qp}
\end{array}\right],$$
where $X_{ij}$ is $m_i \times n_j$ with
$m_1+m_2+\ldots+m_q=n=n_1+n_2+\ldots+n_p.$
In this paper we will require that the diagonal
blocks be square, viz. $p=q$ and
$m_1=n_1,\ m_2=n_2,\ldots ,m_q=n_q.$
Rather than using $m_1, m_2,\ldots,m_q$ to
describe our block structure, we will look at
the sets $\upsilon_i=\{m_{i-1}+1,\ldots m_i\}$.
More formally, we will say the sequence
$\Upsilon =\lla \upsilon_1,\upsilon_2,
\ldots ,\upsilon_{r}\rra$
is an {\it (ordered) partition } of $
\langle n\rangle$, if $\cup_{i=1}^r \upsilon_i=
\langle n\rangle $, and $\upsilon_i \cap \upsilon_j = \emptyset$,
for all $i\neq j$.  We say $\Upsilon$
is an {\it order preserving partition } of $
\langle n\rangle$ if $\Upsilon$
is a partition such that if $i<j$, then $i \in \upsilon_{k_i}$
and $j\in \upsilon_{k_j}$
with $k_i\leq k_j$.
We will say $X$ is {\it block lower triangular}
with respect to the order preserving partition
$\Upsilon=\lla \upsilon_1,\ldots,\upsilon_r  \rra$
if $X_{\upsilon_i\ \upsilon_i^+}=0,$
for every $i\in\langle r\rangle$.
Clearly $X$ may be block lower triangular with respect to several
different partitions.
We call $\Psi =\lla \psi_1,\psi_2,\ldots ,\psi_{p}  \rra$
a {\it refinement}  of an order preserving partition
$\Upsilon$,  if $\Psi$ is also an
order preserving partition of $\langle n\rangle$ and if for every
$i\in\langle p \rangle$ there exists $
j\in\langle r \rangle$ such that
$\psi_i\subseteq \upsilon_j.$ A refinement is said to be {\it  proper}
if $
\Psi\neq\Upsilon.$
The refinement relation on the set of order preserving partitions
of $\langle n \rangle$ defines a lattice such that greatest
lower bound of two partitions
has as its elements intersections of the elements of the
two partitions. The
maximal element is $\lla \{1,2,\ldots,n\} \rra$ , and the minimal element is
$\lla \{1\},\{2\},\ldots,\{ n\}\rra$.  (For the corresponding result for
unordered partitions see
\cite[Lemma 1, p.192]{G} and \cite[Theorem 6, p.7]{B}).
It is easy to see that the greatest lower bound, $\Theta
=\lla \theta_1,\ldots,\theta_q  \rra$, of
the set of partitions for which $X$ is block lower
triangular is also a partition for which $X$ is
block lower triangular, and we will refer to $\Theta$ as
{\it the block lower triangular self--partition of $X$.}
Similarly, the greatest lower bound, $\Phi=
\lla \phi_1,\ldots,\phi_p  \rra$, of the partitions
$\Upsilon =\lla \upsilon_1,\upsilon_2,
\ldots ,\upsilon_{r}\rra$,
for which $X_{\upsilon_i^+\ \upsilon_i}=0$ is referred
to as {\it the block upper triangular self--partition of $X$.}
Thus we say $X$ is {\em lower triangular} if
the block lower triangular self--partition of $X$ is
$\lla \{1\},\{2\}\ldots\{ n\}\rra$.
Similarly, we say
$X$ is {\em upper triangular} if
the block upper triangular self--partition of $X$ is
$\lla \{1\},\{2\}\ldots\{ n\}\rra.$

To illustrate the definitions above, let
$$X=\left[\begin{array}{cccccccc}
*&*&*&0&0&0&0&0\\
0&*&0&0&0&0&0&0\\
0&0&*&0&0&0&0&0\\
0&0&0&*&0&0&*&0\\
0&0&0&0&0&*&0&0 \\
0&0&0&0&0&0&0&0\\
0&0&0&0&*&0&*&0\\
0&0&0&0&0&0&0&*
\end{array}\right],$$
where $*$ denotes a nonzero entry.  Then
viewed as a block lower triangular matrix
$$X=\left[\begin{array}{c|c|c|c|c|c|c|c}
\multicolumn{1}{c}{*}&\multicolumn{1}{c}{*}&*&\multicolumn{1}{c}{0}&
\multicolumn{1}{c}{0}&\multicolumn{1}{c}{0}&\multicolumn{1}{c}{0}&
\multicolumn{1}{c}{0}\\
\multicolumn{1}{c}{0}&\multicolumn{1}{c}{*}&0&\multicolumn{1}{c}{0}&
\multicolumn{1}{c}{0}&\multicolumn{1}{c}{0}&\multicolumn{1}{c}{0}&
\multicolumn{1}{c}{0}\\
\multicolumn{1}{c}{0}&\multicolumn{1}{c}{0}&*&\multicolumn{1}{c}{0}&
\multicolumn{1}{c}{0}&\multicolumn{1}{c}{0}&\multicolumn{1}{c}{0}&
\multicolumn{1}{c}{0}\\
\cline{1-7}
\multicolumn{1}{c}{0}&\multicolumn{1}{c}{0}&0&\multicolumn{1}{c}{*}&
\multicolumn{1}{c}{0}&\multicolumn{1}{c}{0}&*&
\multicolumn{1}{c}{0}\\
\multicolumn{1}{c}{0}&\multicolumn{1}{c}{0}&0&\multicolumn{1}{c}{0}&
\multicolumn{1}{c}{0}&\multicolumn{1}{c}{*}&0&
\multicolumn{1}{c}{0}\\
\multicolumn{1}{c}{0}&\multicolumn{1}{c}{0}&0&\multicolumn{1}{c}{0}&
\multicolumn{1}{c}{0}&\multicolumn{1}{c}{0}&0&
\multicolumn{1}{c}{0}\\
\multicolumn{1}{c}{0}&\multicolumn{1}{c}{0}&0&\multicolumn{1}{c}{0}&
\multicolumn{1}{c}{*}&\multicolumn{1}{c}{0}&*&
\multicolumn{1}{c}{0}\\
\cline{4-8}
\multicolumn{1}{c}{0}&\multicolumn{1}{c}{0}&\multicolumn{1}{c}{0}&
\multicolumn{1}{c}{0}&\multicolumn{1}{c}{0}&\multicolumn{1}{c}{0}&
0&\multicolumn{1}{c}{0}\end{array}\right],$$
and has block lower triangular self--partition
$(\{1,2,3\},\ \{4,5,6,7\},\ \{8\}).$
However, viewed as a block upper triangular matrix
$$X=\left[\begin{array}{c|c|c|c|c|c|c|c}
*&\multicolumn{1}{c}{*}&\multicolumn{1}{c}{*}&\multicolumn{1}{c}{0}&
\multicolumn{1}{c}{0}&\multicolumn{1}{c}{0}&\multicolumn{1}{c}{0}&
\multicolumn{1}{c}{0}\\
\cline{1-2}
0&*&\multicolumn{1}{c}{0}&\multicolumn{1}{c}{0}&
\multicolumn{1}{c}{0}&\multicolumn{1}{c}{0}&\multicolumn{1}{c}{0}&
\multicolumn{1}{c}{0}\\
\cline{2-3}
\multicolumn{1}{c}{0}&0&*&\multicolumn{1}{c}{0}&
\multicolumn{1}{c}{0}&\multicolumn{1}{c}{0}&\multicolumn{1}{c}{0}&
\multicolumn{1}{c}{0}\\
\cline{3-4}
\multicolumn{1}{c}{0}&\multicolumn{1}{c}{0}&0&*&
\multicolumn{1}{c}{0}&\multicolumn{1}{c}{0}&\multicolumn{1}{c}{*}&
\multicolumn{1}{c}{0}\\
\cline{4-7}
\multicolumn{1}{c}{0}&\multicolumn{1}{c}{0}&\multicolumn{1}{c}{0}&0&
\multicolumn{1}{c}{0}&\multicolumn{1}{c}{*}&0&
\multicolumn{1}{c}{0}\\
\multicolumn{1}{c}{0}&\multicolumn{1}{c}{0}&\multicolumn{1}{c}{0}&0&
\multicolumn{1}{c}{0}&\multicolumn{1}{c}{0}&0&
\multicolumn{1}{c}{0}\\
\multicolumn{1}{c}{0}&\multicolumn{1}{c}{0}&\multicolumn{1}{c}{0}&0&
\multicolumn{1}{c}{*}&\multicolumn{1}{c}{0}&*&
\multicolumn{1}{c}{0}\\
\cline{5-8}
\multicolumn{1}{c}{0}&\multicolumn{1}{c}{0}&\multicolumn{1}{c}{0}&
\multicolumn{1}{c}{0}&\multicolumn{1}{c}{0}&\multicolumn{1}{c}{0}&
0&\multicolumn{1}{c}{0}\end{array}\right],$$
and has block upper triangular self--partition
$$(\{1\},\ \{2\},\ \{3\},\ \{4\},\ \{5,6,7\},\ \{8\} ).$$

Let $J_1,\ldots,J_r$ be
subsets of $\langle n \rangle$.
We say that an order preserving partition  $\Psi=\lla
\psi_1,\ldots,\psi_t  \rra$, {\it encompasses }
$J_1,\ldots , J_r$,  if for each $i\in\langle r\rangle
$, there exists $k\in\langle t\rangle$ such that
$J_i\subseteq \psi_k$.
It is easy to see that
the greatest lower
bound of the order preserving partitions which encompass
$J_1,\ldots,J_r$, also encompasses $J_1,\ldots, J_r$ and
we will refer to this as the {\it finest order
preserving partition encompassing} $J_1,\ldots,J_r$.

\bigskip

Next we define subsets
$\at_i$ and $\af_i$ associated with the matrix $A$.
These subsets are defined in terms of the access relationships
in $\G(A)$
to and from the singular classes of $A$, and are convex
in the sense that if $j,k\in \at_i$ (or $\af_i$),
then $l\in \at_i$ ($\af_i$) for
every $j\leq l\leq k$.

\begin{defi} \label{san} {\em Let $A\in\RE ^{nn}$.
Let $S_1,S_2,\ldots S_m $ be
 the singular classes of $A$ ordered so that
 $max(S_i)<max(S_{i+1})$. For each
 $i\in \langle m \rangle$, let
$\mu _i=max(S_i),$ and}
$$\af_i=\{ l\geq\mu_i \  | \ \mbox{there exists $j\geq l$ such that
$j $ is accessed FROM $S_i$ in $\G(A)$} \}. $$
$$\at_i=\{ l\geq\mu_i \  | \ \mbox{there exists $j\geq l$ such that
$j$ has access TO $S_i$ in $\G(A)$} \}, $$
\end{defi}

These subsets turn out to be the key to understanding the
block structure of
a block LU factorization of $A$.

Notice that the $\af_i$ for $A$ correspond to the $\at_i$ for $A^T$ and
vice versa.

\begin{rem}
{\em \label{vct} In \cite{VC}, Varga and Cai show that
an M--matrix admits an LU factorization into M--matrices with $L$ nonsingular
if and only if $\at_i=\{\mu_i\}$ for every $i\in\langle m\rangle$.  }
\end{rem}

 Notice that if $j>i$ and vertices in $S_j$
are accessed from vertices in
$S_i$, then $\af_j\subset \af_i$, so if $i$ is placed in $J$,
then $j$ should be also. Similarly if
vertices in $S_j$ have access to vertices in
$S_i$, then $\at_j\subset \at_i$ and hence if $i$ is placed in $K$, then
$j$ should be also.

In Theorem \ref{sbs} we show that
for each $i\in\langle m\rangle$, either $\af_i$ is encompassed in the
block lower triangular self--partition of $L$ or $\at_i$ is
encompassed in the block upper triangular self--partition of $U$.
This suggests that when factoring
an M--matrix $A$, these $\at_i$ and $\af_i$ should be examined in order to
determine
an optimal factorization of a given type.  There
are several possibilities for the types of factorizations
one might want.  In this paper we highlight four
possibilities.  In Strategy \ref{a} we outline a strategy
for choosing a partition $(J,K)$ so as to minimize the
sizes of the blocks in a block LU factorization of $A$.  In
Strategy \ref{b} we outline a strategy for choosing a permutation matrix
$P$ so that $PAP^T$ has an LU factorization.
In Section 4 we look at the structure of $U$ if $L$ is required to be
nonsingular
(or vice versa).  In Section 5, we look at factoring $A$ as $LBU$ where $L$ is
a nonsingluar lower triangular M--matrix, $U$ is a
nonsingular upper trianglular M--matrix, and $B$ is
block diagonal.

\msn

\begin{exa} {\em
Here we provide an example which
illustrates the definitions introduced above.

Let $$A=\left[ \begin{array}{cccc}
0&-1&0&0\\
0&0&0&-1\\
0&-1&0&0\\
0&0&0&1
\end{array}\right].$$

Then the block lower triangular self--partition
of $A$ is $\lla \{1,2,3,4\}\rra$, and
the block upper triangular self--partition is $\lla\{1\},\{2,3\},\{4\}\rra$.
The singular classes of $A$ are $$S_1=\{1\},\ S_2=\{2\},\
S_3=\{3\},$$ and
$$\at_1=\{1\},\ \at_2=\{2,3\},\ \at_3=\{3\},$$
$$\af_1=\{1,2,3,4\},\ \af_2=\{2,3,4\},\ \af_3=\{3\}.$$
Notice that $\at_3\subseteq \at_2$ and $\af_3\subseteq \af_2\subseteq
\af_1$.  }
\end{exa}

\bigskip
In Section 4, we expand the results of \cite{VC} to
include a larger set of LU factorizations, and
we examine the block structure when no LU factorization exists.
The class associated with an M--matrix
$A$ defined next includes all of the matrices $U$ for which
$A=LU$, with $L,U\in\M$ and $L$ nonsingular.  Notice that
this is actually a larger class since it allows for
some M--matrices $U\in\M$ for which $A=LU$ with $L$ inverse
nonnegative but  not necessarily an M--matrix..

 For any $A\in \RE^{nn}$, we say a matrix  $X$ is
{\it class nonsingular for}
$A$ if for every class $K$ of $A$, $X_{KK}$ is nonsingular
(see  \cite{M}). Note that $K$ need not be a class of $X$.
We write
$$\N = \{\ X\ | \ X\geq 0 \mbox{ and $X$ is class
nonsingular for $A$} \}.$$

\begin{defi}
{\em Let $A$ be an M--matrix and define} $$\U =\{\ BA \ | \
B\in\N, \  BA\in\M \}.$$
$$\UT =\{\ AB \ | \ B\in\N, \  AB\in\M \}.$$
\end{defi}

\section{Block LU Factorizations with (Possibly) Singular $L$ and $U.$}

In \cite{VC}, Varga and Cai consider LU factorizations of M--matrices
where L is nonsingular (see Remark \ref{vct}).  In Example \ref{eg},
the matrix $A$
does not satisfy the conditions stated in \cite{VC} and hence does not
have a
LU factorization into M--matrices with
$L$ is nonsingular.  It does however, have an LU factorization when both
$L$ and $U$ are allowed to be singular.

\begin{exa} \label{eg} {\em Let}
$$A= \left[ \begin{array}{rrr}  0&-1&0\\0&0&0\\0&-1&0 \end{array}
\right]. $$
{\em  By \cite[Theorem 1,(see Remark \ref{vct})]{VC}, $A$  does
not admit a factorization $A=LU$ with $L$ a
nonsingular lower triangular M--matrix and $U$
an upper triangular M--matrix.  By applying the result
to $A^T$ we see that $A$  does not admit a factorization $A=LU$ with
$L$ a lower triangular M--matrix and $U$ a nonsingular
upper triangular M--matrix.  However, consider}
$$L= \left[ \begin{array}{rrr} 1&0&0\\0&0&0\\0&-1&0
\end{array} \right], \ \ \
U= \left[ \begin{array}{rrr}  0&-1&0\\0&1&0\\0&0&0
\end{array} \right]. $$
{\em Then $A=LU$, where $L$ is a singular lower triangular
M--matrix and $U$ is a singular upper
triangular M--matrix.}
\end{exa}

\bigskip

In this section we establish necessary and sufficient conditions
for $A$ to have
an LU factorization into M--matrices,
allowing both $L$ and $U$
to be singular.  This result is stated as Corollary \ref{eq2}.
Enroute to establishing  this result we also characterize the
minimum block structure of block
LU factorizations
when no triangular factorization exists.  These
results appear as Theorem \ref{sbs} and  Theorem \ref{ebs}.

In Strategy \ref{a}, we outline a
strategy one might take to minimize the block sizes
in a block LU factorization of $A$.  In Strategy \ref{b}, we
suggest a permutation $P$ such that
$PAP^T$ admits an LU factorization.

\medskip

We begin with two lemmas and a corollary which we use to prove
the main results in this section (Theorem \ref{sbs}, Theorem \ref{ebs} and
Corollary \ref{eq2}).

\bigskip

First we show that if
$A=LU$ is a factorization of an M--matrix $A$
into M--matrices and $S$ is a singular class of $L$,
then the vertices which are accessed by any vertex of
$S$ in $\G(A)$,
are also accessed by the vertices of $S$ in $\G(L)$.

\begin{lem} \label{ls}
Let $A\in\M$ with factorization $A=LU$
where $ L,U \in\M$.
Let $S$ be any singular class of $L$.  Then
for any $p\in S$ and any $q\in\langle n \rangle$, if
$p$ has access to $q$ in $\G(A)$ then $p$ has access to $q$ in
$\G(L)$.
\end{lem}

\noindent
Proof:

Let $q\in\langle n\rangle$.
Suppose some $p\in S$ has access to $q$ in $\G(A)$ but not in $\G(L)$.
Then by choosing an appropriate section of a path from
$p$ to $q$ in $\G(A)$, we can choose a path $l_1
      \rightarrow l_2 \rightarrow \ldots
 \rightarrow l_g$, in $\G(A)$ so that $$l_1\in S ,$$
$$\mbox{$l_1$ has access to $l_i$ in $\G(L)$ and $l_i \notin S$
for all $i=2,\ldots,g-1$}$$
$$\mbox{$l_1$ does not have access to $l_g$ in $\G(L)$,}$$
$$\mbox{$l_i\neq l_j$, for all $ i\neq j.$}$$

We will establish a contradiction by proving two claims.
The proofs of Claim 1 and Claim 2 rely heavily on the sign
patterns of $A,L,$ and $U$.

\medskip

\begin{description}
\item[Claim 1]  If $r\in S$, $t\notin S$, are such that
$L_{rS}U_{St}<0$, then there exists
$s\in S$ such that $L_{st}\neq 0$.

\medskip

Proof of Claim 1:

\noindent
Since $L_{SS}$ is an irreducible singular M--matrix and
$L_{SS}U_{St}\neq 0,$ by \cite[Theorem 4.16, p. 156]{BP}
it must be a vector which contains both positive and
negative elements.  Hence
there exists $s\in S$ such that $L_{sS}U_{St}>0.$ But then
$0\geq A_{st}=L_{st}U_{tt}+P$, where $P>0$.  Hence
$L_{st}U_{tt}<0,$ and in particular $L_{st}\neq 0.$  This
establishes Claim 1.

\medskip

\item[Claim 2]  For each $j\in\langle g -1\rangle$,

\begin{description}
     \item[(a)] There exists $r\in S\cup\{l_2,\ldots ,l_j\}$

	 with
         $L_{rl_{j+1}}\neq 0$,

    \item[(b)] If $t\notin S\cup\{l_2,\ldots ,l_{j+1}\}$ and
         $U_{l_{j +1}t}\neq 0,$  then there exists

         \noindent
         $r\in S\cup\{l_2,\ldots ,l_j\}$ with $L_{rt}\neq 0,$

\end{description}

\medskip

Proof of Claim 2:

\noindent
We  proceed by induction on $j$.

\begin{description}
\item Let $j=1:$

\begin{description}

\item[(a)]  Since $0>A_{l_1l_2}=L_{l_1l_2}U_{l_2l_2}+L_{l_1S}U_{Sl_2}+P$, where
$P\geq 0,$ either $L_{l_1l_2}U_{l_2l_2}<0$ and
thus $L_{l_1l_2}\neq 0$, or $L_{l_1S}U_{Sl_2}<0$ and thus by
Claim 1 there exists $s\in S$ such that $L_{sl_2}\neq 0$.

\item[(b)] Suppose $U_{l_2t}\neq 0$, for some $t\notin S\cup\{ l_2\}.$
    By (a) there exists $r\in S$ with $L_{rl_2}\neq 0$, hence
$L_{rl_2}U_{l_2t}>0$.  Thus
$0\geq A_{rt}=L_{rt}U_{tt}+
                 L_{rS}U_{St}+P$, where $P>0.$
 			 Hence either
$L_{rt}U_{tt}<0$ and thus $L_{rt}\neq
0$, or $L_{rS}U_{St}<0$ and by
Claim 1, there exists $s\in S$ such that $L_{st}\neq0$.
\end{description}

\item
Let $k<g$.  Assume true for all $j$ with $1\leq j<k$ and show true for $k$.

\begin{description}
\item[(a)]  Since $0>A_{l_kl_{k+1}}=L_{l_kl_{k+1}}U_{l_{k+1}l_{k+1}}+
L_{l_kl_k}U_{l_kl_{k+1}}+P$, where
$P\geq 0,$ either $L_{l_kl_{k+1}}U_{l_{k+1}l_{k+1}}<0$ and
thus $L_{l_kl_{k+1}}\neq 0$, or $L_{l_kl_k}U_{l_kl_{k+1}}<0$ and so by
the induction hypothesis (b) applied with $j=k-1$ and $t=l_{k+1}$,
 there exists  $r\in S\cup\{l_2,\ldots ,l_{k-1}\}$ such that $L_{rl_{k+1}}\neq
 0.$

\item[(b)] Suppose $U_{l_{k+1}t}\neq 0$, for some $t\notin
S\cup \{l_2,\ldots,l_{k+1}\}.$
 By (a) there exists $r\in S\cup \{l_2,\ldots ,l_{k}\}$
such that $L_{rl_{k+1}}\neq 0,$
hence $L_{rl_{k+1}}U_{l_{k+1}t}>0$.  Thus
$0\geq A_{rt}=L_{rt}U_{tt}+
                 L_{rr}U_{rt}+P$, where $P>0.$
Hence either $L_{rt}U_{tt}<0$ and thus $L_{rt}\neq 0$,
or $L_{rr}U_{rt}<0$ and thus $U_{rt}<0$ and so by the
induction hypothesis applied with $j=k-1$, there exists
         $q\in S\cup\{l_2,\ldots ,l_{k-1}\}$ with $L_{qt}\neq 0.$

\end{description}
\end{description}
This establishes Claim 2.
\end{description}

\medskip
By Claim 2, there exists $r\in S\cup\{l_2,\ldots ,l_{g-1}\}$
such
that $L_{rl_g}\neq 0$, but then
there exists $t\in S$ such that $t$ has access
to $l_g$ in $\G(L)$, and since $S$ is a class
of $\G(L)$, it must be that $l_1$ has access to $l_g$
in $\G(L)$.
Contradiction. Hence $p$ must have access to $q$ in
$\G(L)$. {\hfill $\Box$}

\bigskip
By applying Lemma \ref{ls} to $A^T$ we get the following result.

\begin{cor} \label{us}
Let $A\in\M$ with factorization $A=LU$
where $ L,U \in\M$.
Let $S$ be any singular class of $U$.  Then
for any $p\in S$ and any $q\in\langle n \rangle$, if
$p$ is accessed by $q$ in $\G(A)$ then $p$ is accessed by $q$ in
$\G(U)$.
\end{cor}

\noindent
Proof:

Apply Lemma \ref{ls} to $A^T=U^TL^T$. {\hfill $\Box$}

\bigskip
Next we show that every singular class of $A$ has a
corresponding singular class $Q\subseteq S$ in either $L$ or $U$.

\medskip
\begin{lem}\label{lus} Let $A\in\M$ with factorization $A=LU$ where
$L,U\in\M$.  Let $S$ be any singular class of $A$.
Then either $L_{SS}$ is singular or $U_{SS}$ is singular, and
there exists $Q\subseteq S$ such that $Q$ is a singular class
of either $L$ or $U$.
\end{lem}

\noindent
Proof:

$$A_{SS}=L_{SS}U_{SS}+L_{SS'}U_{S'S}=L_{SS}U_{SS} +P,
\mbox{ where } P\geq 0.$$  Thus $L_{SS}U_{SS}=A_{SS}-P$ and
hence is a Z-matrix.  Suppose that $L_{SS}$ and $U_{SS}$ are
both nonsingular.  Then by \cite[$N_{44}$, p. 137]{BP}, $L_{SS}U_{SS}$
is a nonsingular M--matrix, and by \cite[$I_{27}$, p.137]{BP}  there
exists $x \gg 0$ such that $L_{SS}U_{SS}x \gg 0$.  But then
$A_{SS}x=L_{SS}U_{SS}x + Px \gg 0$, which implies that
$A_{SS}$ is a nonsingular M--matrix.  A contradiction.
Hence either $L_{SS}$ is singular or $U_{SS}$ is singular.
Since any proper principal
submatrix of an irreducible M--matrix
is nonsingular, it must be the case that if
$L_{SS}$ (or $U_{SS}$) is singular, then the singular class $Q$
of $L_{SS}$ (or $U_{SS}$) is a singular class of $L$ (or $U$).
{\hfill $\Box$}

\bigskip
The next theorem shows that certain access relationships
must be present either in $L$ or $U$.  In particular, for each $i\in
\langle m\rangle,$ either $F_i$ is encompassed in a block of $L$ or
$T_i$ is encompassed in a block of $U.$

\begin{thm}[Minimality] \label{sbs} Let $A\in\M$ with factorization
$A=LU$ where $L,U\in\M$.  Let $S_i, \at_i, \af_i$, and
$m$ be as in Definition \ref{san}.  Then
there is a  $Q_i \subseteq S_i$
such that $Q_i$ is a singular class of either $L$ or $U.$
If
$Q_i$ is a singular class
of $L$, then $\af_i$ is encompassed in the block lower
triangular self-partition of $L$.
If
$Q_i$ is a singular class
of $U$, then $\at_i$ is encompassed in the block upper
triangular self-partition of $U$.
\end{thm}

\noindent
Proof:

That there exists such a $Q_i$ follows directly from Lemma \ref{lus}.
If $Q_i$ is a singular class of $L$ then by Lemma \ref{ls}, any $j$
which is accessed from $S_i$ (and hence $Q_i$) in $\G(A)$
is also accessed from $Q_i$ in $\G(L),$ thus $\af_i$ is
encompassed in the block lower triangular self-partition of $L.$
Similarly, if $Q_i$ is a singular class of $U$, then
Corollary \ref{us} implies that $\at_i$ is encompassed in the
block upper triangular self--partition of $U.$
 {\hfill $\Box$}

\bigskip
It is natural for one to ask if, once an assignment of the singular classes
between $L$ and $U$ has been chosen, such factorizations can
be achieved.  The following example shows that the
singular classes of $A$ cannot necessarily be divided up
between $L$ and $U$ to suit ones choosing.

\begin{exa}
{\em Let }
$$A=\left[\begin{array}{cc}
0&-1\\
0&0
\end{array}\right].
$$
{\em The singular classes of $A$ are $S_1=\{ 1\}$ and $S_2=\{ 2\}.$
It is easy to verify that there is no factorization
$A=LU$ where $S_1$ is a singular class of $L$ and $S_2$ is a
singular class of $U$. Notice that $F_2\subset F_1$ and hence
there is no combinatorial benefit to having $S_2$ not be a singular
class of $L$ once $S_1$ has been chosen to be a singular class of $L$.}

\end{exa}

In the next theorem we show that once a partition $(J,K)$ of $\langle
m\rangle$ has been
chosen,  a block factorization can be achieved with the block  structure
of $L$ being a refinement of the partition encompassing $F_i,\ i\in J$
and the block structure of $U$ being a refinement of the partition
encompassing $T_i,\ i\in K.$

\begin{thm}[Existence] \label{ebs}
Let $A\in\M$. Let $S_i, \at_i, \af_i$, and $m$ be as in Definition \ref{san}.
Let $(J, K)$ be a partition of $\langle m\rangle$.
Let $\Psi$ be the finest order preserving partition
of $\langle n \rangle$ encompassing
$\af_i,\ i\in J$ and let $\Upsilon$ be the
finest order preserving partition of
$\langle n\rangle$ encompassing $\at_i,\ i\in K$.  Then there exists a
factorization $A=LU$ such that $L,U\in\M$,
the block lower triangular self-partition of $L$ is
a refinement of $\Psi$, and the
block upper triangular self-partition of $U$
is a refinement of $\Upsilon$.
\end{thm}

\noindent
Proof:

\medskip
\noindent

We establish this result by outlining a recursive algorithm which gives us the
desired factorization.

\begin{description}
\item If $a_{11}\neq 0$ then let
$$ N=\{2,\ldots,n\}$$
and set
$$B=A_{NN}-\frac{1}{a_{11}}A_{N1}A_{1N}.$$
Notice that $B$ is an M-matrix and $\G(B)\subseteq \Gamma$ where $\Gamma$
is the subgraph of $\overline{\G(A)}$ induced by the vertices
$2,\ldots,n.$  Moreover, if $Q$ is a singular
class of $B$
(where $B$ is indexed with indices corresponding to $A$) then
$Q\subseteq S$ where $S$ is a singular class of $A$.
Now apply the algorithm to $B$ to obtain $B=\hat{L}\hat{U}$
where $\hat{L}$ and $\hat{U}$ satisfy the theorem.  Set

$$L=\left[ \begin{array}{cc}
1&\begin{array}{ccc}0&\ldots& 0\end{array} \\
\begin{array}{c}\frac{a_{21}}{a_{11}} \\
\frac{a_{31}}{a_{11}} \\
\vdots\\
\frac{a_{n1}}{a_{11}} \end{array} & \hat{L}
\end{array}
\right] \mbox{ and }
U=\left[\begin{array}{cc}
a_{11}&\begin{array}{ccc}a_{12}&\ldots& a_{1n}\end{array} \\
\begin{array}{c}0 \\
0 \\
\vdots\\
0 \end{array} & \hat{U}
\end{array}
\right].
$$
Then $A=LU$ satisfies the theorem.

\item If $a_{11}=0$, then  $1=\mu_1.$
\begin{description}
\item[(i)] If $1\in J$ then set
 $$V=\{ l\ |\  \mu_1\ has\ access\ to\ l\ in \ \G(A)\ \},$$
  $$W=\langle n\rangle \setminus V$$
Notice that $A_{VW}=0,\ V\subseteq \af_1,$ and $max(V)=max(\af_1)$.
Choose a permutation matrix $P$ such that
$$P^{-1}AP=
\left[ \begin{array}{cc}
A_{VV}&0\\
A_{WV}&A_{WW}
\end{array}\right]
=
\left[ \begin{array}{cc}
A_{VV}&0\\
A_{WV}&I
\end{array}\right]
\left[ \begin{array}{cc}
I&0\\
0&A_{WW}
\end{array}\right] ,$$
where elements of $V$ and $W$ are listed in ascending order.
Since $A_{WW}$ is a principal submatrix of $A$ it is also an M--matrix and
we can apply the algorithm to $A_{WW}$ to get $A_{WW}=\hat{L}\hat{U}$ where
$\hat{L}$ and $\hat{U}$ have the desired structure
based on the properties of $A_{WW}$.
Set
$$L=
P\left[ \begin{array}{cc}
A_{VV}&0\\
A_{WV}&\hat{L}
\end{array}\right] P^{-1}, \ \
U=P\left[ \begin{array}{cc}
I&0\\
0&\hat{U}
\end{array}\right]P^{-1}.$$
Then $A=LU.$ We now argue that this factorization satisfies the theorem.
Since applying the permutation similarity only reorders vertices, we see
that $L_{VW}=0$ and $ L_{WW}=\hat{L}$ where $\hat{L}$ has the desired block
structure based on the properties of $A_{WW}.$  We also have that  $V\subset
\af_1,$  thus the block lower triangular self--partition of $L$ is a refinement
of $\Psi.$  By again observing that applying the permutation
similarity with $P$ merely reorders vertices, we see that
$U_{VV}$ is a diagonal matrix, $U_{VW}=0,\
U_{WV}=0$, and $U_{WW}=\hat{U}$ has the desired structure based on the
properties of $A_{WW}.$ Thus the block upper triangular self--partition
of $U$ is a refinement of $\Upsilon.$

\item[(ii)] If $1\notin J$ then $1\in K,$ and we can
apply the algorithm to $A^T$ with $J$ and $K$ interchanged to get
$A^T=\hat{L}\hat{U}.$ Set $L=\hat{U}^T$ and $U=\hat{L}^T.$ Since
transposing a matrix reverses the access relationships, the $\af_i$
for $A$ correspond to the $\at_i$ for $A^T$ and vice versa.  By
the argument presented in $(i)$, the factorization $A=LU$ has
the desired properties.
\end{description}
\end{description}

{\hfill $\Box$}

\medskip
\begin{exa}
{\em Let \label{sense} }
$$A= \left[ \begin{array}{rrr}  0&-1&0\\0&0&0\\0&-1&0 \end{array}
\right]. $$
{\em Let $J=\{1\}$ and $K=\{2,3\}$. Then $(J,K)$ is a partition
of the singular classes of $A$, however since $\af_2\subset \af_1,$
$S_2$ is automatically a singular class of $L$, thus $\at_2$, as the following
factorization shows, need not be encompassed in
the block structure of $U$.
Notice that
$$L= \left[ \begin{array}{rrr} 0&-1&0\\0&0&0\\0&-1&1
\end{array} \right], \ \ \
U= \left[ \begin{array}{rrr}  1&0&0\\0&1&0\\0&0&0
\end{array} \right]. $$
provides an LU factorization which is a proper
refinement of the partition identified by
Theorem \ref{ebs}.  Here $S_2$ is a singular class of $L$
rather than of  $U.$

The partitioning of the singular classes
between $L$ and $U$ of this factorization is
actually $J=\{1,2\},\ K=\{3\}.$  Using this
paritition of the singular classes, Theorem \ref{sbs} and
Theorem \ref{ebs} guarantee that the $L$ and the $U$ listed above
have the smallest possible blocks for this assignment of the
singular classes between $L$ and $U$.}

\end{exa}

\bigskip
Next we state necessary and
sufficient conditions for
an M--matrix $A$ to admit an LU factorization
into M--matrices,
thus extending the results in $\cite{VC}$ to the case where both
$L$ and $U$ are allowed to be singular.

\begin{cor} \label{eq2}
Let $A\in\M$.  Then the following are equivalent:

\begin{description}
\item[(i)]  $A$ admits a factorization $A=LU$, where $L$ is a lower
triangular
M--matrix and $U$ is a upper triangular M--matrix.
\item[(ii)] Let $S_i, \at_i, \af_i$, and $m$ be as in Definition \ref{san}.
Then for every $i\in\langle m \rangle$, either
$\at_i=\{\mu_i\}$ or
$\af_i=\{\mu_i\}$.
\end{description}
\end{cor}

\begin{str} {\em Using \label{a} Theorem \ref{sbs} and Theorem \ref{ebs}
one can strategically choose, for example, a partition which minimizes
the sizes of the blocks in $L$ and $U$.
Let $S_i,\ \mu_i,\ \at_i, \af_i$, and $m$ be as in Definition
\ref{san}.
For each $i\in\langle m\rangle, $ either $\af_i$ has to be encompassed in
block of $L$ or $\at_i$ has to be encompassed in a block of $U$. In each
case, we choose the smallest block between the two, unless one of the
sets ($\at_i$ or $\af_i$) is a subset of an earlier choosen set in
which case it has already been taken care of.

Begin by setting $V=\langle m\rangle.$
\begin{description}
\item[(i)] For $i=min(V)$, if
$| \af_i|<| \at_i|,$  put $i$ into $J$, otherwise put $i\in K$.
Remove $i$ from $V$.  At this time other elements from $V$ may
have automatically been taken care of
(see Example \ref{sense}). Hence if $i$ was placed in $J$, then for
each $j\in V$ such that $\af_j\subseteq \af_i,$ place $j$ in
$J$ and remove $j$ from $V$.
If $i$ was placed in $K$, then for
each $j\in V$ such that $\at_j\subseteq \at_i,$ place $j$ in
$K$ and remove $j$ from $V$.
\item[(ii)] Repeat (i) until $V=\emptyset$.
\item[(iii)]  Apply the algorithm provided in the proof of Theorem
\ref{ebs} with the partition $(J,K).$
\end{description}
 }
\end{str}

\begin{str} {\em Another \label{b} strategy one might employ is to choose a
permutation matrix $P$ such that $PAP^T$ satisfies Corollary
\ref{ebs}, and hence has an LU factorization.  There are several ways
one might do this.  For example, using
$\mu_i$ and $m$ as in Definition
\ref{san}, one could choose a permutation matrix $P$ which
corresponds to the permutation which moves $\mu_1,
\ldots \mu_m$ to positions $n-m+1, \ldots n$
and reorders them
(if necessary) so that $\mu_i$ is placed after
$\mu_j$ whenever $\mu_i$ has access to $\mu_j$ in $\G(A).$ All other
indices should be shifted up by the appropriate amount.
Then $\af_i=\{\mu_i\}$ so by Corollary \ref{eq2}, a triangular
LU factorization exists with $J=\langle m\rangle$ and $K=\emptyset.$
}

\end{str}

We conclude this section with two examples on which we illustrate
the strategies suggested by the theorems in this section.

\begin{exa} \label{eg2}
{\em Let }
$$A=\left[\begin{array}{ccccccc}
1& -1& 0& 0& -1& 0&0\\
-1& 1& 0& 0& -2& 0&0\\
0 &0&2& -2& 0& 0&0\\
0&0& -2&2&0&0&0\\
0&0&0&0&1&0&0 \\
0&0&-1&-1&0&0&0\\
0&0&0&0&0&-1&1
\end{array}\right].$$
{\em Let $S_i,\ \mu_i,\ \at_i,\ \af_i$ be as in Definition \ref{san}.
Then}
$$S_1=\{ 1,2\},\ S_2=\{ 3,4\},\ S_3=\{6\},$$
$$\mu_1=2,\ \mu_2=4,\ \mu_3=6,$$
 $$ \at_1=\{2\} ,\ \at_2=\{ 4,5,6,7\},\ \at_3=\{6,7\},$$
 $$ \af_1=\{2,3,4,5\} ,\ \af_2=\{4\},\ \af_3=\{6\}. $$
{\em  Thus by \cite[Theorem 1]{VC} (see Remark \ref{vct})
$A$ does not admit an LU factorization into M--matrices
with $L$ nonsingular.  Similarly by \cite[Theorem 1]{VC} applied to
$A^T$,
$A$ does not admit an LU factorization into M--matrices
with $U$ nonsingular.  Using Strategy \ref{a} we choose $J=\{2,3\}$ and
$K=\{1\}.$  The algorithm outlined in the proof of Theorem \ref{ebs}
now gives us an LU factorization with }
$$L= \left[\begin{array}{ccccccc}
1& 0& 0& 0& 0& 0&0\\
-1& 1& 0& 0& 0& 0&0\\
0 &0&2& 0& 0& 0&0\\
0&0&-2&0&0&0&0\\
0&0&0&0&1&0&0 \\
0&0&-1&-2&0&0&0\\
0&0&0&0&0&-1&1
\end{array}\right],\  U=\left[\begin{array}{ccccccc}
1& -1& 0& 0& -1& 0&0\\
0& 0& 0& 0& -3& 0&0\\
0 &0&1& -1& 0& 0&0\\
0&0& 0&1&0&0&0\\
0&0&0&0&1&0&0 \\
0&0&0&0&0&1&0\\
0&0&0&0&0&0&1
\end{array}\right].$$
\end{exa}

\begin{exa}
{\em Let \label{eg3}}
$$A=\left[\begin{array}{cccccccc}
1&-1&0&0&0&0&0&0\\
0&0&-1&0&0&0&0&0\\
0&0&1&0&-1&0&0&0\\
0&0&0&1&0&0&0&0\\
0&0&-1&0&1&-1&0&0\\
0&0&0&0&0&0&0&0 \\
0&0&0&0&0&0&0&-1\\
-1&0&0&0&0&0&0&0
\end{array}\right].$$
{\em Let $S_i,\ \mu_i,\ \at_i,\ \af_i$ be as in Definition \ref{san}.
Then }
$$S_1=\{ 2\},\ S_2=\{ 3,5\},\ S_3=\{6\},\ S_4=\{7\},\ S_5=\{8\}, $$
$$\mu_1=2,\ \mu_2=5,\ \mu_3=6,\ \mu_4=7,\ \mu_5=8,$$
 $$ \af_1=\{2,3,4,5,6\} ,\ \af_2=\{5,6\},\ \af_3=\{6\},\
 \af_4=\{7,8\},\ \af_5=\{8\}, $$
 $$ \at_1=\{1,2,3,4,5,6,7,8\},\at_2=\{5,6,7,8\}, \at_3=\{6,7,8\},
\at_4=\{7\}, \at_5=\{8\}. $$
{\em Using strategy \ref{a} one should set $J=\{1,2,3\},$ and
$K=\{4,5\}.$   Using the algorithm outlined in Theorem \ref{ebs}
we get }
$$L=\left[\begin{array}{cccccccc}
1&0&0&0&0&0&0&0\\
0&0&-1&0&0&0&0&0\\
0&0&1&0&-1&0&0&0\\
0&0&0&1&0&0&0&0 \\
0&0&-1&0&1&-1&0&0\\
0&0&0&0&0&0&0&0 \\
0&0&0&0&0&0&1&0\\
-1&-1&0&0&0&0&0&1
\end{array}\right],\
U=\left[\begin{array}{cccccccc}
1&-1&0&0&0&0&0&0\\
0&1&0&0&0&0&0&0\\
0&0&1&0&0&0&0&0\\
0&0&0&1&0&0&0&0\\
0&0&0&0&1&0&0&0 \\
0&0&0&0&0&1&0&0\\
0&0&0&0&0&0&0&-1 \\
0&0&0&0&0&0&0&0
\end{array}\right].$$
{\em This factorization minimizes the sizes of the
blocks in $L$ and $U$.}

{\em Using strategy \ref{b}, we can choose
a permutation which moves $2,5,6,7,8$ to the end of our list of indices
and then reorders them as $6,5,2,8,7$.  Thus we take }
$$P=\left[\begin{array}{cccccccc}
1&0&0&0&0&0&0&0\\
0&0&1&0&0&0&0&0\\
0&0&0&1&0&0&0&0\\
0&0&0&0&0&1&0&0\\
0&0&0&0&1&0&0&0 \\
0&1&0&0&0&0&0&0\\
0&0&0&0&0&0&0&1 \\
0&0&0&0&0&0&1&0
\end{array}\right].$$ {\em Then }
$$PAP^{-1}=\left[\begin{array}{cccccccc}
1&0&0&0&0&-1&0&0\\
0&1&0&0&-1&0&0&0\\
0&0&1&0&0&0&0&0\\
0&0&0&0&0&0&0&0\\
0&-1&0&-1&1&0&0&0\\
0&-1&0&0&0&0&0&0\\
-1&0&0&0&0&0&0&0\\
0&0&0&0&0&0&-1&0\\
\end{array}\right],$$ {\em which has an LU factorization with}
$$L=\left[\begin{array}{cccccccc}
1&0&0&0&0&0&0&0\\
0&1&0&0&0&0&0&0\\
0&0&1&0&0&0&0&0\\
0&0&0&0&0&0&0&0\\
0&-1&0&-1&0&0&0&0\\
0&-1&0&0&-1&0&0&0\\
-1&0&0&0&0&-1&0&0\\
0&0&0&0&0&0&-1&0\\
\end{array}\right],\
U=\left[\begin{array}{cccccccc}
1&0&0&0&0&-1&0&0\\
0&1&0&0&-1&0&0&0\\
0&0&1&0&0&0&0&0\\
0&0&0&1&0&0&0&0\\
0&0&0&0&1&0&0&0\\
0&0&0&0&0&1&0&0\\
0&0&0&0&0&0&1&0\\
0&0&0&0&0&0&0&1\\
\end{array}\right].$$

\end{exa}


\section{Block LU Factorizations with Nonsingular $L$}

In this section we consider the case where an
M--matrix
$A$ is factored into block
triangular matrices with
$L$ nonsingular, such that $L^{-1}$ is nonnegative
and class nonsingular for $A$, and $U\in\M$.  Notice
that if $L$ is a nonsingular M--matrix, then $L$ satisfies
the given conditions, thus we are considering a larger set
of LU factorizations than were considered in \cite{VC}.
We would like $L$ and $U$ to be as
close to lower and upper
triangular respectively as possible.  We begin by
showing that certain access relationships
must be present in  $\G(U)$.

\begin{thm}[Minimality] \label{dtf} Let $A\in\M$ and
$U\in\U$.  Let $S_1,S_2,\ldots S_m $ be
 the singular classes of $A$.  Then  for each
 $i\in \langle m \rangle$ and $j\in\langle n\rangle$,
if  $j$ has access to
$S_i$ in $\G(A)$, then
$j$ has access to some vertex $q\in S_i$
in $\G(U)$.
\end{thm}

Proof: Assume $j$ has access to
$S_i$ in $\G(A)$. Since $A_{S_iS_i}$ is
irreducible, $j$ has access in $\G(A)$ to every vertex of $S_i$.
Let $Q_i$ be any
final class of $U_{S_iS_i}.$
For any $q\in Q_i$, we see that $j$ has access to $q$ in
$\G(A)$.   Choose
$B\in\N$ such that $BA=U$. Applying
\cite[Theorem 3.7]{M} we see that
$j$ has access to $q$ in $\G(U)$.
{\hfill $\Box$}

\bigskip

Next we show that there is a  $U\in\U$, for
which the only nonzeros below the diagonal
of $U$ correspond to
the required access relationships described in
Theorem \ref{dtf}. To
optimize on the placement of these zeros, for
each $S_i$, we
would like the vertex $q$ in Theorem \ref{dtf} to be
$\mu_i$.  The desired  $U$
can be found by performing Gaussian elimination until
a zero diagonal element is encountered.   When a
zero diagonal element is encountered, that column should be
skipped over, and the process continued. Thus we obtain a matrix $U$
which is upper triangular except for spurs below the diagonal corresponding to
the columns $\mu_i$. We thus obtain a $U$ whose column envelope is small,
where {\em column envelope} is defined similarly to the row envelope in
\cite[p. 708]{GP}.

Notice that using this algorithm the
$L$ which is produced is
a nonsingular lower triangular M--matrix.

\begin{thm}[Existence]  \label{tt} Let $A\in\M$.  Let
 $S_1,S_2,\ldots S_m $ be
 the singular classes of $A$,  and let
$\mu _i=max(S_i).$
Set
 $$\chi=\{(j,\mu_i) \ | \ i\in\langle m \rangle, \  j>\mu_i
\mbox{ and $j$ has access to $\mu_i$ in } G(A)\}.$$
Then $A$ admits a
factorization $A=LU$, where $L\in\M$ is lower
triangular and nonsingular, and $U\in\M$ is such that
for all $j>k$, $u_{jk}= 0$ whenever $(j,k)\notin\chi$,
and $u_{jj}=0$ if and only if $j=\mu_i $ for some $i\in\langle m\rangle$.
\end{thm}

Proof:
Let $A^\varepsilon$ be the matrix obtained
from $A$ by adding $\varepsilon$ to the
$(\mu_i,\mu_i)-th$ element in $A$, for
each $i\in\langle m\rangle$.
Then for any $\varepsilon>0$ and
any $i\in\langle m \rangle$, $A^\varepsilon_{S_iS_i}$ is nonsingular by
\cite[Theorem 4.16(2) and
Theorem 2.7]{BP}.  Hence
 $A^\varepsilon $
is a nonsingular M--matrix and thus
all the diagonal
elements of $A^\varepsilon $ are positive.
For each $i\not\in
\{ \mu_1,\ldots,\mu_m\}$, and
for each $j>i$, use multiplication by elementary matrices
on the left
to add appropriate multiples of row $i$ to
row $j$ so as to reduce $A^\varepsilon $ to a matrix for which
the $(j,i)-th$ element is $0$.
Since
the matrix formed at
each step of this process is
also a  nonsingular M--matrix,
we will not encounter a nonzero diagonal element as we
proceed in this fashion.

Notice that since
row $\mu_k, k\in\langle m \rangle$, is not used as
a pivot row, the
elementary matrices will be independent of
$\varepsilon$. By proceeding in this manner
we can produce a matrix $U^\varepsilon$ such that
$u^\varepsilon_{ji}=0$ whenever
$i\not\in \{ \mu_1,\ldots,\mu_m\}$ and $j>i$.  Moreover,
$\G(U^\varepsilon)\subseteq\overline{\G(A)}$.

The off diagonal elements of $U^\varepsilon $ remain nonpositive, hence
$U^\varepsilon$ is a Z--matrix.
Since the leading principal minors of $A^\varepsilon $ are positive, and
adding a multiple of one row to another does not change the
determinant, all the leading principal minors of $U^\varepsilon$ are
positive. Hence  $U^\varepsilon$ is a nonsingular M--matrix.
Since $\varepsilon$ only appears on the diagonal of $U^\varepsilon$,
$U=\lim_{\varepsilon\rightarrow 0 }U^\epsilon$  is an
M--matrix such that
$u_{jk}= 0$ whenever $(j,k)\notin\chi$.
It is
easy to see that if we take $L$ to be the inverse of the
product of the elementary matrices used above (in the appropriate order)
then $L$
is a nonsingular lower triangular M--matrix and
$A=LU$.
Since all of the diagonal elements of $A$ where
used as pivots, except for the $(\mu_i,\mu_i)-th$,
they must be nonzero both in $A$ and in $U$.
Since $A_{S_iS_i}$ is
singular, it must be that
$U_{S_iS_i}$ is singular.
By the construction of
$U$, $U_{S_i S_i}$
is an upper triangular matrix with
nonzeros on the diagonal except possibly for
$u_{\mu_i\mu_i}$, hence
it must be that
$u_{\mu_i\mu_i}=0$.
Thus $U$ has the desired pattern.
{\hfill $\Box$}

\begin{exa}
{\em Let}
$$A=\left[\begin{array}{cccccccc}
     1&    -1  &   0 &    0  &   0   &  0  &   0  &   0\\
    -1 &    1  &   0 &    0  &   0  &  -1   &  0 &    0\\
    -1 &    0   &  2   & -2   &  0  &   0   & -1   &  0\\
     0 &   -1   & -1   &  1 &    0 &    0   &  -2   & -1\\
     0  &   0   &  0  &  -1 &    1  &   0  &   0 &    0\\
     0   &  0 &    0  &   0   &  0  &   0   &  0   &  0\\
     0   &  0 &    0  &   0  &   0    &-1   &  1   &  0\\
     0   &  0 &    0 &    0&     0&    -1 &   -1 &    1
\end{array}\right]$$
{\em
Then $\mu_1=2, \ \mu_2=4,\ \mu_3=6$ and
$\chi=\{(3,2),(4,2),(5,2),(5,4),(7,6),(8,6)\}$
Using our algorithm we get}
$$L=\left[ \begin{array}{cccccccc}
1 &   0 &   0  &     0  & 0 & 0   &   0&0\\
-1  & 1  &  0  &    0   & 0 & 0   &   0&0\\
-1  & 0  &  1  &    0  &  0 &  0  &    0&0\\
0 &0  & -\frac{1}{2}&  1 &  0&  0 &0&0\\
0  & 0  &  0   &    0  &1  &   0    &     0&0\\
0  & 0 &   0    &   0  & 0   & 1    &     0&0\\
0  & 0  & 0    &    0  & 0   &  0  &  1&0\\
0  & 0 &  0    &    0  &0    &  0 &  -1&1
\end{array}\right] \mbox{, \em and }
U =\left[ \begin{array}{c|c|c|c|c|c|cc}
\multicolumn{1}{c}{1}&  \multicolumn{1}{c}{-1}  &\multicolumn{1}{c}{0}&
\multicolumn{1}{c}{0} & \multicolumn{1}{c}{0} & \multicolumn{1}{c}{0}
& 0&0\\
\multicolumn{1}{c}{0}& \multicolumn{1}{c}{0}& \multicolumn{1}{c}{0} &
\multicolumn{1}{c}{0} &\multicolumn{1}{c}{0} &\multicolumn{1}{c}{-1}& 0&0\\
\cline{2-2}
 0 &  -1   & \multicolumn{1}{c}{2} & \multicolumn{1}{c}{-2}
 &\multicolumn{1}{c}{0}&
 \multicolumn{1}{c}{0}& -1&0\\
0 &  -\frac{3}{2}    & \multicolumn{1}{c}{0}&\multicolumn{1}{c}{0}&
\multicolumn{1}{c}{0}&\multicolumn{1}{c}{0}&
\multicolumn{1}{c}{-\frac{5}{2}}&-1\\
\cline{4-4}
 0& 0&0& -1& \multicolumn{1}{c}{1}  &\multicolumn{1}{c}{0} &    0&0\\
 \cline{2-2}\cline{4-4}
\multicolumn{1}{c}{0}&\multicolumn{1}{c}{0}   &\multicolumn{1}{c}{0}&
\multicolumn{1}{c}{0}  & \multicolumn{1}{c}{0}  &\multicolumn{1}{c}{0}& 0&0\\
\cline{6-6}
 \multicolumn{1}{c}{0}   &\multicolumn{1}{c}{0}  & \multicolumn{1}{c}{0}
 & \multicolumn{1}{c}{0}&     0 &  -1    &1&0\\
\multicolumn{1}{c}{0}   & \multicolumn{1}{c}{0}   & \multicolumn{1}{c}{0}
&\multicolumn{1}{c}{0} &    0&   -2     & 0&1\\
\cline{6-6}
\end{array}\right].$$         {\em The boxes indicate the
positions of potential nonzero elements below the diagonal of
$U$.  Notice
that even though $(5,2)\in\chi,$ $u_{5,2} = 0$.  Theorem
\ref{tt} only asserts that the elements in $\chi$ are
potentially nonzero.  }
\end{exa}

\begin{rem} \label{zero}
{\em Theorem \ref{dtf} gives a lower bound on the number of
nonzeros which
must occur below the diagonal of $U$. Let
$S_1, \ldots S_m $  be the singular classes of $A$.
Let $\mu_i=max(S_i)$, $i\in\langle m\rangle$.
Let $$R_i=\{ \ j\in\langle n\rangle\ |\ j>\mu_i\ and\ j\ has\
access\ to\ \mu_i\ in\ \G(A)\}.$$
Let $$R=\cup_{i=1}^n R_i.$$  Then $|R|\leq$ the number of nonzeros below the
diagonal of $U$.  Moreover, Theorem \ref{tt}
shows that there is a factorization of $A$ into
a nonsingular lower triangular M--matrix $L$ and
a block upper triangular M--matrix $U$, for which the
 number of nonzeros below
the diagonal of $U$ $\leq \sum_{i=1}^m |R_i|$.
The next example shows that if
the $R_i$ have nonempty pairwise intersection, then the
number of nonzeros corresponding to this duplication
is not identified by our theorems.  }
\end{rem}

\begin{exa} {\em

Consider
$$A=\left[\begin{array}{ccc}
0&0&0\\
-1&0&0\\
0&-1&1\\
\end{array}\right],\ \ \
B=\left[\begin{array}{ccc}
0&0&0\\
-1&0&0\\
-1&-1&1\\
\end{array}\right].$$
Both $A$ and $B$ have singular classes $S_1=\{1\},\ S_2=\{2\}$, and
in both $\G(A)$ and $\G(B)$, vertices $2$ and $3$ access $S_1$, and
vertex $3$ accesses $S_2$.
Using the notation of Remark \ref{zero},  $R_1=\{2,3\},\
R_2=\{3\},$ and $R=\{2,3\}.$  By Theorem \ref{dtf},
there must be at least 2 nonzero  subdiagonal elements
in $U$ for any block LU factorization of $A$ or $B$ (where
$L$ and $U$ satisfy the conditions of
Theorem \ref{dtf}), and by Theorem \ref{tt},
there is an LU factorization of $A$ and an LU factorization of
$B$, each
with at most $3$ subdiagonal elements in the $U$ and with the
$L$ a nonsingular lower triangular M--matrix.
Notice that if we
factor $A=IA$, then $U=A$ has $2$ nonzeros below the diagonal.
However it can
be shown that for any factorization of $B=LU$ (where $L$ and $U$ satisfy
the conditions of Theorem \ref{dtf}), $U$ will have
at least $3$ nonzeros
below the diagonal.
}
\end{exa}

\medskip
Our next theorem  characterizes the M--matrices $A$ such that
$A=LU$ where $L,U\in\M$ with $L$ nonsingular and lower triangular and
$U$ upper triangular.
The equivalence of parts (ii) and (iii) of this theorem
appear in \cite{VC}.

\begin{thm} \label{eq}
Let $A\in\M$.  Then the following are equivalent:

\begin{description}
\item[(i)]  There exists $U\in\U$ such that $U$ is upper triangular.
\item[(ii)]  $A$ admits a factorization $A=LU$, where $L$ is a
nonsingular lower triangular
M--matrix and $U$ is an upper triangular M--matrix.
\item[(iii)]  Let $m,T_i,$ and $\mu_i$ be
as in Definition \ref{san}.  Then
$\at_i=\{\mu_i\}, \ \forall i\in\langle m\rangle$.
\end{description}
\end{thm}

\noindent
Proof:

Follows from Theorem \ref{dtf} and Theorem \ref{tt}.
{\hfill $\Box$}

\bigskip

\section{Block LBU Factorizations}

Using Theorem \ref{tt} we can factor
$A=LBU$ into the product of a nonsingular lower triangular M--matrix
$L$
with
a block diagonal M--matrix $B$ and a nonsingular
upper triangular M--matrix $U$.
Since $L$ and $U$
are nonsingular lower and upper triangular M--matrices,
standard methods can be used to solve the parts of the
system in which they are involved.
In general, $B$ will be a sparse matrix since its only
nonzero off--diagonal entries occur in rows and columns corresponding
to the ends of singular classes.  In this specialised case it
may be possible to use sparse matrix techniques to
solve the part of the system involving $B$.

\begin{thm} \label{lbu} Let $A\in\M$.  Let $$\chi=\{ \ (j,j)\
| \ j\in\langle n \rangle\setminus \{ \ \mu_i\  | \  i\in\langle
m \rangle \}\ \} $$
$$
\cup
\{\ (j,\mu_i) \  | \  i\in\langle m \rangle, \ j>\mu_i \mbox{  and $j$
has access to $\mu_i$ in $\G(A)$}\}$$
$$\cup
\{ \ (\mu_i,j) \  | \   i\in\langle m \rangle, \ j>\mu_i \mbox{  and $j$
is accessed by $\mu_i$ in $\G(A)$}\}.$$  Then $A$ admits a
factorization $A=LBU$ where $L$ is a nonsingular lower
triangular M--matrix, $U$ is a nonsingular upper triangular
M--matrix, and $B$ is an M--matrix  such that
$b_{jk}=0$ whenever $(j,k)\notin\chi$.  Moreover, if $L$ and $U$ are
nonsingular M--matrices, this is the finest block structure that $B$
can have (i.e it is impossible to subdivide the blocks of $B$
and obtain a block diagonal matrix).

\end{thm}

\noindent
Proof:

\medskip
\noindent
Use Theorem \ref{tt} to factor $A=LV$ where $L$ is a nonsingular
M--matrix and $V$ is an M--matrix
such that
for all $j>k$, $v_{jk}= 0$ whenever $(j,k)\not\in
\{(j,\mu_i) \  | \  i\in\langle m \rangle, \ j>\mu_i \mbox{  and $j$
has access to $\mu_i$ in $\G(A)$}\} $. Notice that the singular
classes of $V$ are just $\mu_1,\ldots,\mu_m$, and $\G(V)\subseteq
\overline{\G(A)}$.
Using the algorithm
outlined in the proof of Theorem \ref{tt}, we can factor $V^T=XY$ where
$X$ is an nonsingular lower triangular M--matrix, and $Y$ is
an M--matrix such that for all $j>k$, $y_{jk}=0$ whenever $(j,k)\not\in
\{(j,\mu_i) \  | \  i\in\langle m \rangle, \ j>\mu_i \mbox{  and $j$
has access to $\mu_i$ in $\G(A)$}\} $. Since
rows $\mu_1,\ldots,\mu_m$ are not used as pivot rows,
and all other rows of $V^T$
have zeros to the right of the diagonal, $Y$ will also
satisfy the property that for all $j>k$, $y_{kj}= 0$ whenever
$(j,k)\not\in
\{(j,\mu_i) \  | \  i\in\langle m \rangle, \ j>\mu_i \mbox{  and $j$
has access to $\mu_i$ in $\G(A)$}\}. $
Let $B=Y^T$ and $U=X^T$.  Then $A=LBU$ is as claimed.  By Theorem
\ref{dtf}, the block structure exhibited by $B$ is the finest
possible with $L$ and $U$ nonsingular M--matrices.
{\hfill $\Box$}

\begin{exa}
{\em Let \label{eg4}}
$$A=\left[\begin{array}{cccccccc}
1&-1&0&0&0&0&0&0\\
0&0&-1&0&0&0&0&0\\
0&0&1&0&-1&0&0&0\\
0&0&0&1&0&0&0&0\\
0&0&-1&0&1&-1&0&0\\
0&0&0&0&0&0&0&0 \\
0&0&0&0&0&0&0&-1\\
-1&0&0&0&0&0&0&0
\end{array}\right],$$ {\em as in Example \ref{eg3}
Then }
$$\mu_1=2,\ \mu_2=5,\ \mu_3=6,\ \mu_4=7,\ \mu_5=8,$$
{\em and}
$$\chi=\{(2,3),(2,4),(2,5),(2,6),(7,2),(8,2),(5,6),(7,5),(8,5)
(7,6),(8,6),(7,8)\}.$$

{\em Using the algorithm outlined above Theorem \ref{tt}
we get $A=LV$.   Then using the algorithm to factor $V^T=U^TB^T$
we get }
$$L=\left[\begin{array}{cccccccc}
1&0&0&0&0&0&0&0\\
0&1&0&0&0&0&0&0\\
0&0&1&0&0&0&0&0\\
0&0&0&1&0&0&0&0 \\
0&0&-1&0&1&0&0&0\\
0&0&0&0&0&1&0&0 \\
0&0&0&0&0&0&1&0\\
-1&0&0&0&0&0&0&1
\end{array}\right],\
B=\left[\begin{array}{c|c|cc|c|c|c|c|}
\multicolumn{1}{c}{1}&\multicolumn{1}{c}{0}&0&\multicolumn{1}{c}{0}
&\multicolumn{1}{c}{0}&\multicolumn{1}{c}{0}&\multicolumn{1}{c}{0}
&\multicolumn{1}{c}{0}
\\
\cline{3-6}
\multicolumn{1}{c}{0}&0 &-1&
\multicolumn{1}{c}{0}
&\multicolumn{1}{c}{-1}&0&
\multicolumn{1}{c}{0}
&\multicolumn{1}{c}{0}
\\
\cline{3-6}
\multicolumn{1}{c}{0}
&\multicolumn{1}{c}{0} &1&\multicolumn{1}{c}{0}&\multicolumn{1}{c}{0}&
\multicolumn{1}{c}{0}&
\multicolumn{1}{c}{0}
&\multicolumn{1}{c}{0}
\\
\multicolumn{1}{c}{0}
&\multicolumn{1}{c}{0} &0&\multicolumn{1}{c}{1}&\multicolumn{1}{c}{0}
&\multicolumn{1}{c}{0}&
\multicolumn{1}{c}{0}
&\multicolumn{1}{c}{0}
\\
\cline{6-6}
\multicolumn{1}{c}{0}
&\multicolumn{1}{c}{0}&0&\multicolumn{1}{c}{0}&0&-1&
\multicolumn{1}{c}{0}
&\multicolumn{1}{c}{0}
\\
\cline{6-6}
\multicolumn{1}{c}{0}
&\multicolumn{1}{c}{0} &0&\multicolumn{1}{c}{0}&\multicolumn{1}{c}{0} &
\multicolumn{1}{c}{0}
&\multicolumn{1}{c}{0}&\multicolumn{1}{c}{0}
\\
\cline{2-2}\cline{5-6}\cline{8-8}
0&0&0&0&0&0&0&-1\\
\cline{8-8}
0&-1&0&0&0&
0&\multicolumn{1}{c}{0}
&\multicolumn{1}{c}{0}
\\
\cline{2-2}\cline{5-6}
\end{array}\right],\ $$
$$U=\left[\begin{array}{cccccccc}
1&-1&0&0&0&0&0&0\\
0&1&0&0&0&0&0&0\\
0&0&1&0&-1&0&0&0\\
0&0&0&1&0&0&0&0\\
0&0&0&0&1&0&0&0 \\
0&0&0&0&0&1&0&0\\
0&0&0&0&0&0&1&0 \\
0&0&0&0&0&0&0&1
\end{array}\right].$$
{\em The boxes in $B$ indicate the only possible positions for nonzero
off  diagonal entries.}
\end{exa}

\bigskip
\bigskip

\centerline{{\bf Acknowledgments}}

We would like to thank Michael Tsatsomeros
and an anonymous referee
for helpful comments and
suggestions.

\end{document}